\let\ORIlabel\label
\let\ORIrefstepcounter\refstepcounter
\AddToHook{package/hyperref/before}{
    \let\label\ORIlabel
    \let\refstepcounter\ORIrefstepcounter
}
\documentclass[final,onefignum,onetabnum]{siamart220329}
\usepackage{amsfonts,amsmath,mathrsfs,bm} 
\usepackage{amssymb}
\usepackage[notref,notcite]{showkeys}
\usepackage{hyperref}
\usepackage{graphicx}
\usepackage{wrapfig}
\usepackage{subfig}
\usepackage{float}
\usepackage{bbm}
\usepackage{algorithmic}
\usepackage{gensymb,color,soul}
\usepackage{tikz}
\usepackage{tikz-3dplot}
\usepackage{booktabs}
\usepackage{mathtools}
\usepackage{comment}
\mathtoolsset{showonlyrefs}

\newcommand{\N}{\mathbb{N}}
\newcommand{\R}{\mathbb{R}}

\newtheorem{remark}[theorem]{Remark}

\newcommand{\vertiii}[1]{{\left\vert\kern-0.25ex\left\vert\kern-0.25ex\left\vert #1 
\right\vert\kern-0.25ex\right\vert\kern-0.25ex\right\vert}}

\makeatletter
\newcommand*\bigcdot{\mathpalette\bigcdot@{.5}}
\newcommand*\bigcdot@[2]{\mathbin{\vcenter{\hbox{\scalebox{#2}{$\m@th#1\bullet$}}}}}
\makeatother

\title{Projections for handling uncertainties and enabling domain truncation in diffuse optical tomography}

\author{
A.~Hakula\footnotemark[1]
\and P.~Hirvi\footnotemark[1]
\and N.~Hyv\"onen\footnotemark[1] \footnotemark[3]
\and A.~Jääskeläinen\footnotemark[1]
\and V.~Kolehmainen\footnotemark[2]
}

\begin{document}
\maketitle

\renewcommand{\thefootnote}{\fnsymbol{footnote}}
\footnotetext[1]{Aalto University, Department of Mathematics and Systems Analysis, P.O.~Box 11100, FI-00076 Aalto, Finland (aada.hakula@aalto.fi, pauliina.hirvi@aalto.fi, nuutti.hyvonen@aalto.fi, altti.jaaskelainen@aalto.fi). The work of PH, NH and AJ was supported by the Research Council of Finland (Flagship of Advanced Mathematics for Sensing Imaging and Modelling grant 359181). The work of AH was part of the Ministry of Education and Culture’s Doctoral Education Pilot under Decision No. VN/3137/2024-OKM-6 (Doctoral Education Pilot for Mathematics of Sensing, Imaging and Modelling).}
\footnotetext[2]{University of Eastern Finland, Department of Technical Physics, Kuopio Campus, P.O. Box 1627, FI-70211 Kuopio, Finland (ville.kolehmainen@uef.fi). The work of VK was supported by the Research Council of Finland (Flagship of Advanced Mathematics for Sensing Imaging and Modelling grant 358944).}
\footnotetext[3]{Corresponding author.}

\begin{abstract} 
This paper presents a projection-based technique to mitigate the impact of modeling errors related to domain truncation, changes in the optode coupling coefficients, and misspecified optical parameters of different tissue types in diffuse optical tomography. The approach considers the primary Jacobian matrix of the forward map in the image reconstruction scheme, linking the primary unknown,~i.e.,~the per-voxel absorption coefficient changes in the region of interest, to the optode measurements, as well as the nuisance Jacobians that do the same for the auxiliary unknown parameters of secondary interest. To mitigate mismodeled coupling coefficients or domain truncation, the method projects the linearized forward model defined by the primary Jacobian onto the orthogonal complement of the range of a nuisance Jacobian, or onto the orthogonal complement of the span of a number of first left singular vectors for the nuisance Jacobian that has been weighted to account for prior information on the measurement setup. In the case of a misspecified baseline optical parameter for some tissue type, the nullspace of the utilized orthogonal projection is defined to be the span of first left singular vectors for a (weighted) difference of two Jacobian matrices evaluated at two different levels for the considered tissue-wise optical parameter. The reconstruction is formed by applying Bayesian inversion with Gaussian prior and noise models to the projected linearized equation. We evaluate the method on simulated brain activity data obtained via Monte Carlo simulations of the radiative transfer equation in a voxelized head anatomy for a neonate with combined gestational and chronological age of 41.7 weeks.
\end{abstract}

\renewcommand{\thefootnote}{\arabic{footnote}}



\pagestyle{myheadings}
\thispagestyle{plain}
\markboth{A.~HAKULA, P.~HIRVI, N.~HYV\"ONEN,  A.~J\"A\"ASKEL\"AINEN AND V. KOLEHMAINEN}{PROJECTIONS FOR MARGINALIZATION IN DOT}

\section{Introduction}
\label{sec:introduction}
This work considers reducing the effect of uncertainties in the forward model of functional diffuse optical tomography (DOT) which maps the absorption coefficient changes to changes in the measurements~\cite{arridge1999optical,Nissila_DOI2005}. The basic idea is to project the linearized forward model onto the orthogonal complement of a subspace that is expected to be affected the most by the uncertainties and to subsequently apply Bayesian inversion to the projected equation. A similar approach has previously been considered in electrical impedance tomography (EIT) \cite{Jaaskelainen25,Jaaskelainen26} and analyzed theoretically in \cite{Calvetti2025}.

In this study, we focus on DOT as a developing functional neuroimaging technique for mapping hemodynamic activity on the cerebral cortex. It can be used to estimate changes in oxygenated (HbO$_\mathrm{2}$), deoxygenated (HbR), and total (HbT) hemoglobin concentration by analyzing variations in visible red to near‑infrared light intensity after it traverses tissue~\cite{Nissila_DOI2005,bluestone2001three}. Synaptic activity influences arterial diameter, perfusion, and volume via neurovascular coupling, while shifts in oxygen metabolism alter the balance of hemoglobin types~\cite{hillman2007depth}. Together, these vascular and metabolic processes produce the observed hemodynamic signal.

In this study, we consider DOT in the frequency-domain (FD) and use stochastic Monte Carlo (MC) simulations in a voxel-based head model to predict measurements of log-amplitude and phase delay given the optical properties and the optode (source and detector) configuration. The forward solution, captured through the trajectories and pathlength statistics of detected photons, also produces the Jacobian matrix that linearly relates the measured changes in log-amplitude and phase to perturbations in the medium’s voxel-wise absorption coefficients~\cite{Hirvi2026}. This approximate linear model is then inverted to estimate the absorption coefficient changes in the {\em region of interest} (ROI) from the corresponding difference data. The reconstructed absorption changes at two wavelengths can subsequently be converted to changes in HbO$_\mathrm{2}$ and HbR for physiological interpretation~\cite{bluestone2001three,gibson2005recent}. Alternatively, a single wavelength of 798\,nm is sufficient to estimate the changes in HbT~\cite{cope1991development,maria2022}. This wavelength is selected for the current study to make it physically meaningful to inspect only one reconstruction of each target.

We consider three types of uncertainties in the forward model of DOT: (i) changes in the optode coupling coefficients that model the amplitude losses and phase shifts specific to individual source and detector channels~\cite{Schweiger2007, Mozumder2013}, (ii) absorption coefficient changes outside the selected ROI, and (iii) uncertainty in the tissue-specific baseline optical parameters that are used in the computation of the Jacobian matrix. The ultimate goal is to enable reconstruction of the absorption changes in the ROI only (e.g., brain cortex), even if there are also changes elsewhere in the tissues,~i.e.,~in the {\em region of non-interest} (RONI), or in the coupling coefficients, or if the baseline optical properties of some tissue types were misspecified when forming the primary Jacobian matrix. Variations in the coupling coefficients may be related to,~e.g.,~air gaps between optodes and the skin, movement of the optic fibers, movement of hair under the optodes, or blood flow in superficial blood vessels~\cite{Schweiger2007,Stott2003,Nissila_DOI2005}. The aim to only reconstruct the absorption changes in the ROI, even if there may also be activity in the RONI (e.g., in the extracerebral tissue~\cite{AE_heiskala2012} or untargeted regions in the brain) that affects the measurements, is to speed up the online phase of the reconstructions process by allowing truncation of the computational domain and to potentially provide extra stability,~cf.~\cite{Jaaskelainen26}. This can be highly beneficial, for example, when processing long time series of data. Furthermore, in a real-world experiment, it may be of interest to select the ROI based on the performed task or provided stimulus, for example, limiting to the occipital lobe when imaging responses to visual stimuli~\cite{liao2012high}. Finally, the motivation to tackle the case of misspecified baseline optical parameters stems from the wide range of reference values found in literature \cite{arridge1999optical,farina2015invivo,maria2022,mozumder2024diffuse}, implying that the primary Jacobian may often have to be computed with imprecise knowledge of the baseline optical properties. The computational load would also be reduced if measurements at multiple wavelengths could be handled with Jacobians simulated for one set of baseline optical parameters.  

To address uncertainties in the coupling coefficients, we form the {\em nuisance Jacobian matrix} of the optode measurements with respect to the coupling coefficients, project the linear equation relating the absorption changes to the difference measurements onto the orthogonal complement of the range of the nuisance Jacobian, and compute the reconstruction of the absorption changes by applying Bayesian inversion to the projected equation,~cf.~\cite{Jaaskelainen25}. As the number of coupling coefficients and the number of measurements are, respectively, proportional to the number of optodes, and its square, the projected equation is nontrivial, and as evidenced by our numerical examples, it allows efficient artifact reduction leading to reconstructions that are similar to the ideal case when the optode coupling coefficients are exactly known.

When considering domain truncation, the nuisance parameters are the per-voxel absorption coefficient changes in the RONI, and thus the corresponding Jacobian’s range is much higher-dimensional, often spanning nearly the entire image space. A full projection onto its orthogonal complement would therefore eliminate almost all measurement information. To address this, we employ a partial projection: we select only a subspace of the nuisance Jacobian’s range that is expected to have the strongest impact on the data and project the to-be-inverted ROI-only forward model onto its orthogonal complement. Such a subspace can be formed by computing a number of left singular vectors for a weighted nuisance Jacobian, with the weight incorporating a prior model about the absorption changes within the RONI,~cf.~\cite{Jaaskelainen26}.

The projection for reducing the effect of a misspecified optical baseline parameter for some tissue type is not based on a Jacobian matrix with respect to that parameter. Instead, one computes a weighted difference of Jacobian matrices with respect to the discretized absorption in the ROI for any two literature-based values of the baseline parameter in the considered tissue type and forms the nullspace of the employed orthogonal projection from a number of left singular vectors for this difference matrix. According to our knowledge, such an approach to tackling uncertainties in baseline parameters has not previously been tested in connection to any imaging modality.

Approximate marginalization of uncertainties in the forward model of DOT has previously been addressed via the approximation error (AE) method \cite{kaipio2005book}; see,~e.g.,~\cite{Arridge_2006,Kolehmainen2009,Mozumder2013,Tarvainen2010,Tarvainen2009,AE_heiskala2012}. In the AE approach, the impact of uncertainties is modeled as an additive Gaussian approximation error term in the measurement model. Before any data are collected, the second-order statistics of this error are estimated by simulating random draws from the priors for the free parameters and propagating them through the forward model. The resulting approximation error distribution is then combined with standard Gaussian measurement noise in Bayesian inversion. Recent approaches for handling modeling errors also include deep learning methods~\cite{mozumder2024diffuse} and principal component analysis~\cite{Hakula2026}, which (loosely speaking) are also based on extensive sampling. In contrast, our projection method avoids both additional prior sampling and the corresponding forward solves, which is a principal advantage over the AE and other approaches.

We test and verify our projection approach using FD difference data simulated in a neonate's voxel-based head model with two or five local increases in the absorption coefficient. All difference data and Jacobian matrices are simulated with Monte Carlo eXtreme (MCX)~\cite{fang2009monte, MCX}, which is an open-source software enabling massively parallel computation on a Graphics Processing Unit (GPU). The head model is based on the database presented in \cite{collins2021construction}. In addition to considering the three types of uncertainties in the forward model of DOT separately, we also briefly investigate the possibility of using projections to simultaneously reduce the combined effects of multiple types of uncertainties.

This text is organized as follows. Section~\ref{sec:DOT} introduces DOT and outlines how MCX can be used for simulating DOT measurements and computing related Jacobian matrices. The construction of projections for alleviating the effect of model uncertainties is considered in Section~\ref{sec:projections}, whereas Section~\ref{sec:background} introduces the Bayesian inversion framework, with and without projections. The experimental setup is described in Section~\ref{sec:experimental}. The results of our numerical experiments are documented in Section~\ref{sec:numerics} and the conclusions are drawn in Section~\ref{sec:conclusion}.

\section{Diffuse optical tomography}
\label{sec:DOT}
 
The forward problem in the context of optical imaging usually refers to solving the three-dimensional distribution of radiance or fluence in the domain, given the source locations, as well as the geometry of the domain and distribution of optical properties inside it. Detector locations can then be used to estimate different measurements on the boundary~\cite{arridge1999optical,Nissila_DOI2005}. While Maxwell’s equations are the most fundamental model, the radiative transfer equation (RTE), or its diffusion approximation, is typically adopted for practical computations in turbid media~\cite{arridge1999optical, Nissila_DOI2005}. Since solving the RTE in complex three-dimensional domains is challenging, forward solutions are often obtained with MC methods, which are considered the gold standard for implementing the RTE~\cite{bluestone2001three,fang2009monte}. The accuracy of MC depends on the number of detected photons and can suffer from a low signal-to-noise ratio (SNR), though modern GPUs have made such MC simulations far more practical~\cite{fang2009monte}.

The imaged domain, which is a segmented neonatal head model in this work, is modeled as a bounded regular enough domain $\Omega \subset \R^3$ with a connected complement. For MC, the domain is typically discretized into voxels, which can slightly alter the shape of $\partial \Omega$ and the tissue segmentation. Sources and detectors are identified with connected subsets of $\partial \Omega$. In MCX, the optical properties of $\Omega$ are given by four positive functions from
\[
L^\infty_+(\Omega) = \left\{ v \in L^\infty(\Omega) \; | \; {\rm ess} \!\inf v > 0 \right\}\, ,
\]
which define the absorption coefficient $\mu_a$, the scattering coefficient $\mu_s$, the anisotropy coefficient $g$, and the refractive index $\nu$~\cite{arridge1999optical}. More information on the experimental setup considered in this study is provided in Section~\ref{sec:experimental}.

When there is no activation in the brain causing changes in the optical properties, all coefficients are modeled as piecewise constant over this segmentation of $\Omega$. In functional neuroimaging, it is typically assumed that the absorption coefficient is most altered by hemodynamic activity, compared to the other coefficients~\cite{Nissila_DOI2005,maria2022}. Consequently, we only consider absorption changes in this work.

\subsection{Monte Carlo simulations}
\label{sec:MC}

The MC forward model solves the RTE stochastically by simulating the trajectories of a high number of photon packets originating from an active source on the domain boundary. Since MCX implements precise ray  tracing~\cite{fang2022mcxcloud,FangMMC} and the microscopic Beer--Lambert law~\cite{Sassaroli2012mBBL} to model photon transport, the simulated trajectories are independent of the baseline absorption coefficients and can be reused~\cite{YaoReplay,boas2002three}. At each scattering, the new step length and direction are sampled from the local voxel’s scattering $\mu_s$ and anisotropy coefficients~\cite{fang2009monte, wang1995mcml, boas2002three}. If $\mu_s$ changes mid-step at a tissue boundary, the remaining length is rescaled accordingly~\cite{boas2002three}. At the boundary of $\Omega$, photons either exit or reflect according to Fresnel’s law. Photons exiting on a detector patch are recorded together with their detector index, random number generator seeds required to repeat the same trajectories, and tissue-wise path lengths, providing all data needed for subsequent analysis, i.e., estimating the measurements and simulating the Jacobians~\cite{Hirvi_2023}.

Let $f$ be the intensity-modulation frequency of the FD device, and $l_{p,j}$ be the length of the intersection of voxel $j$ and the trajectory of a photon packet $p$ exiting through a detector. The complex weight of the detected photon packet allows the representation~\cite{heiskala2009accurate,heiskala2007optical,leino2019valomc,Hirvi_2023,kangasniemi2024stochastic}
\begin{equation}
\label{eq:weight}
w_{p} = \exp\!\Big({-\sum_{j} \Big( \mu_{a,j} - \mathrm{i}\ \frac{2\pi f \nu_{j}}{c_0} \Big) \, l_{p,j}}\Big),
\end{equation}
where the initial weight of the packet is assumed to be $1$, $\mu_{a,j}$ is the absorption coefficient in voxel $j$ and $c_0$ is the speed of light in vacuum. Here $\lvert w_{p} \rvert$ approximates the probability that a photon is not absorbed on the prescribed path, and the expected value for this probability at the detector over all simulated photon packets gives an estimate for the relative intensity measured by the source--detector pair in question compared to other pairs. The frequency-dependent term introduces effectively an extra absorption term due to the modulation of input light. Denoting the total time-of-flight for the photon packet $p$ by $t_{p}$, the real and imaginary parts of the complex relative intensity can be separated as~\cite{heiskala2009accurate,heiskala2007optical,Hirvi_2023,Hirvi2026}
\begin{align}
\label{eq:X_Y}
X = \dfrac{1}{N} \sum_{p} \lvert w_{p} \rvert \cos{(2\pi f t_{p})} \quad {\rm and} \quad
Y = \dfrac{1}{N} \sum_{p} \lvert w_{p} \rvert \sin{(2\pi f t_{p})},
\end{align}
where $N$ is the total number of launched photon packets. Typically, FD instruments measure the log-amplitude $\ln(A)$ and phase shift $\varphi$ of the detected photon density wave, which can subsequently be computed via~\cite{Hirvi_2023} 
\begin{equation}
\label{eq:amp_phase}
\ln(A) = \ln \! \big( \sqrt{X^2+Y^2} \big)
\quad {\rm and} \quad
\varphi = \mathrm{atan2} \, \Big(\frac{Y}{X}\Big),
\end{equation}
respectively.

A complete set of measurements for the considered anatomy is simulated by considering in turns all sources and recording the corresponding intensities at the sensors satisfying certain geometric restrictions. Hence, a single set of measurements can be stored,~e.g.,~as a vector $z \in \R^{2 m}$, where $m$ is the number of source--detector pairs used for collecting the data. That is, each source--detector pair produces a single complex number $X + {\rm i}Y$, which is converted into two real measurements via \eqref{eq:amp_phase}. In the following, we assume that the first $m$ elements in $z$ carry the log-amplitude measurements and the last $m$ elements correspond to the phase measurements.

\subsection{Linearized reconstruction model}
\label{sec:linear}
Solving the linearized reconstruction problem of DOT requires the Jacobian of the measurements with respect to voxel-wise absorption. This Jacobian is evaluated at the baseline optical parameters and depends on the tissue segmentation and the measurement configuration, enabling recovery of absorption changes at each time point compared to a reference measurement. 
We refer to \cite[Section~2.2.1]{Hirvi_2023} for more information on how such Jacobian can be computed in MCX; see also \cite{Hirvi2026}, where Jacobians with respect to the scattering coefficient are considered as well.  

Let us denote by $x_{\rm total} \in \R^{n_{\rm total}}$ the voxel-wise absorption change in the whole {\em field of view} (FOV),~i.e.,~the region in $\Omega$ that is expected to have potential to influence the measurements and will be defined more precisely in Section~\ref{sec:FOV}. In anticipation of only being interested in reconstructing the absorption change in the ROI, we decompose $x_{\rm total} = [x^\top, \tilde{x}^\top]^\top$, where $x \in \R^{n}$ corresponds to the voxels in the ROI and $\tilde{x}  \in \R^{\tilde{n}}$, with $n + \tilde{n} = n_{\rm total}$, to those in the RONI that is defined as the complement of ROI relative to the FOV. Analogously, we decompose the Jacobian matrix of the measurements with respect to the voxel-wise absorption change as $J_{\rm total} = [J, \tilde{J}]$, where $J \in \R^{2m \times n}$ and $\tilde{J} \in \R^{2m \times \tilde{n}}$.

Let $y = z - z_0 \in \R^{2m}$ be difference measurement with $z_0$ being the baseline measurement and $z$ corresponding to a measurement at one time point during a functional brain imaging experiment, where the aim is to recover absorption changes caused by hemodynamic variations. Typically, the relation between the absorption changes $x \in \R^{n}$ and the difference data $y \in \R^{2m}$ is approximated by the linear model 
\begin{equation} \label{eq:linear_inverse_problem0}
    y = J x_{\rm total} + e = J x + \tilde{J} \tilde{x} + e
\end{equation}
where $e \in \R^{2m}$ models the measurement noise. In the subsequent sections, we ignore the term $\tilde{J} \tilde{x}$ and simply aim to reconstruct $x$ from 
\begin{equation}
\label{eq:linear_inverse_problem}
y = J x + e.
\end{equation}
Take note that we can return to the unbiased linearized model \eqref{eq:linear_inverse_problem0} by simply defining ROI to be the whole FOV.

\subsection{Coupling coefficients}
\label{sec:ccs}
In addition to the absorption change in the FOV, the difference measurements recorded for source--detector pairs may be influenced by alterations in the coupling coefficients characterizing the optodes~\cite{Schweiger2007,Stott2003,Nissila_DOI2005,mozumder2013optode}. The purpose of this section is to present the related mathematical model and explain how the Jacobian matrix of the measurements with respect to the logarithms of the coupling coefficients can be formed.

Each source and detector have their own coupling coefficient, denoted by $c^{(s)}$ and~$c^{(d)}$, respectively.
These can be written as
\begin{equation}
    c^{(s)} = A^{(s)} \exp ( \mathrm{i} \, \varphi^{(s)} ), \quad c^{(d)} = A^{(d)} \exp ( \mathrm{i} \, \varphi^{(d)} ),
\end{equation}
where $A^{(\cdot)}$ and $\varphi^{(\cdot)}$ are, respectively, the amplitude and the phase of the considered coupling coefficient. We say that optodes are ideally coupled when $A^{(s)}~=~A^{(d)}~=~1$ and $\varphi^{(s)} = \varphi^{(d)} = 0$, i.e., there is no unknown amplitude loss nor phase delay happening in either one of the optodes \cite{Schweiger2007}.

Let us write the uncontaminated absolute measurement corresponding to the source $i$ and detector $j$ as
\begin{equation} \label{eq:measurement_without_coupling_error}
    w_{ij} = A_{ij} \exp ( \mathrm{i} \, \varphi_{ij} ),
\end{equation}
where $A_{ij}$ denotes the amplitude of the measurement and $\varphi_{ij}$ its phase.
Accounting for the coupling coefficients of the source $i$ and the detector $j$, the measurement \eqref{eq:measurement_without_coupling_error} becomes
\begin{equation}
    \tilde{w}_{ij} = c_{i}^{(s)} c_{j}^{(d)} w_{ij} = A_{ij} A_{i}^{(s)} A_{j}^{(d)} \exp \big( \mathrm{i} \, (\varphi_{ij} + \varphi_{i}^{(s)} + \varphi_{j}^{(d)}) \big),
\end{equation}
and taking the logarithm leads to
\begin{equation}
\label{eq:linear_coupling}
    \ln {\tilde{w}_{ij} } = 
    \ln{ A_{ij} } + \ln{ A_{i}^{(s)} } +  \ln{ A_{j}^{(d)} } + \mathrm{i} \big (\varphi_{ij} + \varphi_{i}^{(s)} + \varphi_{j}^{(d)} \big),
\end{equation}
revealing a linear dependence. 

According to \eqref{eq:linear_coupling}, the log-amplitude measurement (respectively, the phase measurement) for each source--detector pair depends linearly on the log-amplitudes (respectively, the phases) of the coupling coefficients at the considered source and detector. Hence, the relation between the difference measurements $y$ and the changes in the log-amplitudes and the phases of the coupling coefficients can be exactly given  with the help of a matrix $J_c \in \R^{2 m \times 2l}$, where $l \in \N$ denotes the total number of sources and detectors. Moreover, $J_c$ is independent of all optical parameters and only carries ones and zeros as its elements, with the corresponding sparsity structure only depending on how the measurements for source--detector pairs (indexed by $i$ and $j$ above) are arranged into the measurement vector $y$.

In what follows, we call $J_c$ the coupling Jacobian and note that it can be formed explicitly without any need for MC simulations. The same applies to the projection onto the orthogonal complement of the range $\mathcal{R}(J_c)$ of $J_c$ considered in the following sections.

\section{Construction of projections}
\label{sec:projections}
In this section, we describe how to construct projections that can mitigate effects of different types of nuisance parameters on the measured data. The rationale behind such projections is that we wish to consider those parts of the measurement data that are not caused by changes in the values of nuisance parameters and utilize projections for reduction of the nuisance parameter related components in the data.  

For the linearized problem, the space of possible changes in the measured data caused by perturbations in certain model parameters is given by the range of the Jacobian matrix computed with respect to those parameters. Therefore, the effects of nuisance parameters can be removed by projecting the data and the forward operator onto the orthogonal complement of the range of the associated Jacobian matrix. The projection can lead to partial or nearly perfect removal of the nuisance effects depending on the Jacobian matrices for the primary and nuisance parameters. The approach can give almost complete removal of the nuisance parameter errors for a low-dimensional nuisance parameter such as the optode coupling coefficients. For high-dimensional parameters, such as the absorption changes outside the ROI (e.g., ROI being part of the brain cortex and RONI being the rest of the head), full projection onto the orthogonal complement of the range of the nuisance Jacobian could eliminate most, if not all of the data. Hence, in such cases we choose a subspace of the range of the nuisance Jacobian for constructing the projection.

\subsection{Low-dimensional parameter}
\label{sec:low}

For a nuisance parameter of a low enough dimension, a matrix that projects away its influence can be computed as follows. Let $J_c \in \R^{2m \times K}$ be the Jacobian matrix of the measurements, computed with respect to the considered parameter and evaluated at some initial guess for all involved parameters. Here, $K \in \N$ is the dimension of the nuisance parameter, for which we assume that $K \ll 2m$. In our experiments, the low-dimensional nuisance parameter corresponds to the coupling coefficients, in which case $K \propto \sqrt{m}$. Then, a projection matrix $P \in \R^{2m \times 2m}$, projecting onto the orthogonal complement of the range of $J_c$, can be formed using the equation
\begin{equation}
\label{eq:proj}
    P = \mathrm{I} - J_c (J_c^{\rm \top}\! J_c)^{-1} J_c^{\rm \top}
\end{equation}
if $J_c$ has full range, which applies to case of coupling coefficients. As the number of employed source--detector pairs $m$ is not typically very high and we assume that $K$ is even smaller, constructing $P$ and operating with it is computationally cheap.

\subsection{High-dimensional parameter}
\label{sec:high}

In the case of a high-dimensional nuisance parameter, the range of the projection must be chosen more carefully. Let us now denote the Jacobian with respect to the nuisance parameter by $\tilde{J}$. In our experiments, this high-dimensional nuisance parameter is the voxel-wise absorption in the RONI, in which case $\tilde{J}$ is defined as described in Section~\ref{sec:linear}. The problem with projecting away too much of the data is solved by choosing a suitable subspace of the range of~$\tilde{J}$, and constructing a matrix that projects onto the orthogonal complement of that subspace.

In order to select the subspace such that we project away as much of the contribution of the nuisance parameter as possible while keeping the dimension of the subspace as low as possible, we construct a basis for the subspace out of directions we expect to be most affected by the nuisance parameter. A simple method for achieving this is choosing the basis vectors to be such orthonormal $\tilde{v}$ that sequentially maximize the norm $\| \tilde{J} \tilde{v} \|_2$. This corresponds to forming the basis out of a chosen number of first left singular vectors of $\tilde{J}$. The same basis can also be constructed by choosing a set of eigenvectors for $\tilde{J} \tilde{J}^{\rm \top}$ corresponding to the largest eigenvalues. Collecting the basis vectors as columns of a matrix $\tilde{V}$, the corresponding projection matrix can be constructed by substituting $\tilde{V}$ for $J_c$ in equation \eqref{eq:proj}.

We may also include additional prior information in the construction of the projection if it is,~e.g.,~known that some alterations in the nuisance parameters are {\em a priori} more likely than others. This can be achieved by introducing a positive semidefinite weighting matrix~$\tilde{A}$, which incorporates available prior information, and forming the basis for the nullspace of the projection out of eigenvectors of $\tilde{J} \tilde{A} \tilde{J}^{\rm \top}$ corresponding to the largest eigenvalues. A possible choice for $\tilde{A}$ is the covariance matrix of an underlying zero-mean prior distribution for the change in the considered nuisance parameter. In this case, the eigenvectors give the directions of maximal variance for $\tilde{J} \tilde{X}$, with $\tilde{X}$ denoting the (randomized) nuisance parameter in question.

\subsection{Inaccurate forward operator}
\label{sec:uncertain}
Assume that the true linearized measurement model is described by $J_\Delta$, whereas $J = J_0$ corresponds to a model based on, say, false prior information on the absorption level of some tissue type that is off by $-\Delta \in \R$. We write the accurate linearized model as
\begin{equation}
\label{eq:uncertain}
J_\Delta x = (J_\Delta - J) x  + J x  = y
\end{equation}
and define
\begin{equation}
\label{eq:differenceJ}
\hat{J}_\Delta  = J_\Delta - J_0.
\end{equation}
Now, one can, in principle, introduce a projection $P_\Delta$ onto the orthogonal complement of $\mathcal{R}(\hat{J}_{\Delta})$ and apply it to \eqref{eq:uncertain} to deduce
\begin{equation}
\label{eq:Puncertain}
P_\Delta J x  =   P_\Delta J_\Delta x - P_\Delta \hat{J}_{\Delta} x = P_\Delta J_\Delta x = P_\Delta y,
\end{equation}
where we have eliminated the effect of the misspecified model parameter on the Jacobian but it has reappeared as a multiplication by $P_\Delta$.

When trying to solve \eqref{eq:Puncertain}, there are two obvious problems: (i)~Being able to form $P_\Delta \in \R$ seems to require knowing $\Delta \in \R$, which would also enable computing $J_\Delta$ and thus considering the accurate linearized forward model \eqref{eq:uncertain} to begin with. (ii)~Even if we knew $P_\Delta$, its nullspace may be so large that the projected equation \eqref{eq:Puncertain} is either trivial or essentially uninformative on $x$. 

Our attempt to solve the first problem is to replace $P_{\Delta}$ by $P_{\delta}$ with a user-defined $\delta \in \R$ that does not, in general, equal $\Delta$. The underlying hope is that although $\hat{J}_{\Delta}$ and $\hat{J}_{\delta}$ are definitely different, the same does not necessarily apply in the same extent to their ranges --- or rather their left singular vectors --- that are the objects defining $P_{\Delta}$ and $P_{\delta}$,~cf.~\cite{Jaaskelainen25}. The second problem is addressed in the same way as in Section~\ref{sec:high}: instead of projecting onto $\mathcal{R}(\hat{J}_{\delta})^{\perp}$, we project onto the orthogonal complement of a subspace spanned by eigenvectors of $\hat{J}_{\delta} A \hat{J}_{\delta}^\top$ corresponding to a number of largest eigenvalues. Here, the symmetric positive semidefinite matrix $A$ incorporates prior information on the unknown variable of primary interest (i.e.,~the absorption change in the ROI).

If there is uncertainty about multiple baseline parameters indexed, say, by $i=1, \dots, q$, we propose to form the basis for the nullspace of the employed orthogonal projection from eigenvectors of
\begin{equation}
\label{eq:sum_matrix}
\sum_{i=1}^q \hat{J}_{i, \delta_i} A \hat{J}_{i, \delta_i}^\top,
\end{equation}
corresponding to the largest eigenvalues. Here, $\hat{J}_{i, \delta_i}$ has been formed as in \eqref{eq:differenceJ}, but with $\Delta$ replaced by $\delta_i$ that is the user-defined perturbation in the $i$th baseline parameter (only). Note that the relative sizes of $\delta_i$, $i = 1, \dots, q$, affect the eigensystem of the sum matrix in \eqref{eq:sum_matrix}, and thus they should be chosen with care based on prior information on the accuracy of the initial estimates for the considered parameters.

\subsection{Combining multiple projections}
\label{sec:combo}
To explain how to simultaneously project away several sources of uncertainties, let $J_c$ and $\tilde{V}$ be as in Sections~\ref{sec:low} and \ref{sec:high}, respectively, and assume the matrix $V$ carries as its columns some number of eigenvectors of the matrix \eqref{eq:sum_matrix} considered in Section~\ref{sec:uncertain}, corresponding to its largest eigenvalues. A projection that approximately projects away the directions affected the most by the types of uncertainties considered in Sections~\ref{sec:low}--\ref{sec:uncertain} is obtained by replacing $J_c$ in \eqref{eq:proj} by
\begin{equation}
\label{eq:manyP}
\mathcal{J} = \big[ J_c, \tilde{V}, V \big],
\end{equation}
assuming its columns are linearly independent --- as they are in all our numerical experiments. If this were not the case, the orthogonal projection could be constructed,~e.g.,~by resorting to a QR decomposition of $\mathcal{J}$.

\section{Bayesian reconstruction algorithm}
\label{sec:background}
In the Bayesian framework, all the uncertainly known variables are modeled as random variables and \eqref{eq:linear_inverse_problem} is rewritten as
\begin{equation}
\label{eq:Bayesian_system}
    Y = JX + E,
\end{equation}
where $Y$, $X$ and $E$ are randomized versions of $y$, $x$ and $e$, respectively. We model the absorption change $X$ and the measurement noise $E$ as mutually independent zero-mean Gaussians with probability distributions
\begin{equation} \label{eq:priors}
        X \sim \mathcal{N}(0,\Gamma_{x}), \quad
        E \sim \mathcal{N}(0,\Gamma_{e}), 
\end{equation}
where $\Gamma_x \in \R^{n \times n}$ and $\Gamma_e \in \R^{2 m \times 2 m}$ are symmetric positive definite covariance matrices. The posterior probability density of the unknown variable $x$ can be written as
\begin{align}
\label{eq:posterior}
    \pi(x \mid y) &\propto \pi_{Y \mid X}(y \mid x)\pi_{X}(x) = \pi_{E}(y - Jx) \pi_{X}(x) \nonumber \\[1mm]
    &\propto \exp \! \Big( -\frac{1}{2} \big( \| y - Jx \|^2_{\Gamma_{e}^{-1}} + \| x \|^2_{\Gamma_{x}^{-1}} \big) \Big),
\end{align}
where the excluded normalization constants do not depend on $x$ and we use the standard notation $\| v \|_{B}^2 = v^\top \! B v$ for a symmetric positive definite matrix $B$. It can be shown that \eqref{eq:posterior} corresponds to a Gaussian density whose (posterior) mean and covariance are~\cite{kaipio2005book}
\begin{align}
\label{eq:mean}
    \widehat{x}_{\mathrm{post}} &= \Gamma_{x}J^\top (J \Gamma_{x} J^\top + \Gamma_{e})^{-1} y, \\[1mm]
  \label{eq:cov}  \Gamma_{\mathrm{post}} &= \Gamma_{x} - \Gamma_{x} J^\top (J \Gamma_{x} J^\top + \Gamma_{e})^{-1} J \Gamma_{x},
\end{align}
which provide the reconstruction and a spread estimator for the absorption change in the ROI when no projections are utilized in the reconstruction process.

Consider then an orthogonal projection $P \in \R^{2m \times 2m}$ and the projected equation
\begin{equation}
\label{eq:proj_eq}
    PY = PJX + PE.
\end{equation}
Note that $PE$ follows a degenerate Gaussian distribution $\mathcal{N}(0, P \Gamma_{e} P^\top)$ on $\mathcal{R}(P)$. In particular, the measurement model \eqref{eq:proj_eq} leads to the posterior density
\begin{align}
\label{eq:Pposterior}
    \pi(x \mid P y) &\propto \pi_{PE}(Py - PJx) \pi_{X}(x) \nonumber \\[1mm]
    &\propto \exp \! \Big( -\frac{1}{2} \big( \| Py - PJx \|_{\Gamma_{e}^{-1}}^2 + \| x \|_{\Gamma_{x}^{-1}}^2 \big) \Big),
\end{align}
where the last step follows from the fact that $Py - PJx$ belongs to the range of $P$, the pseudoinverse of $P$ satisfies $P^{\dagger} = P$, and $P^2 = P$. By completing the square in the argument of the exponential function on the right-hand side of \eqref{eq:Pposterior}, one straightforwardly concludes that $\pi(x \mid P y)$ is a Gaussian density with mean and covariance
\begin{align}
\widehat{x}_{\mathrm{post},P} &=  \Gamma_{\mathrm{post},P}  J^\top \! P^\top \! \Gamma_e^{-1} P y, \\[1mm]
    \Gamma_{\mathrm{post},P} &= \big(\Gamma_{x}^{-1} + J^\top \!P^\top \Gamma_e^{-1} P J\big)^{-1}.
\end{align}
Applying the Woodbury matrix identity to $\Gamma_{\mathrm{post},P}$, it follows via some algebra that 
\begin{align}
\label{eq:Pmean}
    \widehat{x}_{\mathrm{post}, P} &= \Gamma_{x} (PJ)^\top \! \big(PJ \Gamma_{x} (PJ)^\top + \Gamma_{e} \big)^{-1}  P y, \\[1mm]
    \label{eq:Pcov}
    \Gamma_{\mathrm{post},P} &= \Gamma_{x} - \Gamma_{x} (PJ)^\top \! \big(PJ \Gamma_{x} (PJ)^\top + \Gamma_{e} \big)^{-1}P J \Gamma_{x},
\end{align}
which are computationally more efficient formulations as the matrix to be inverted has the dimension $2m \ll n$. These formulas provide the reconstruction and the associated spread estimator when an orthogonal projection is employed in the reconstruction process for mitigating uncertainly known nuisance parameters.

\begin{remark}
It follows straightforwardly from \eqref{eq:cov} and \eqref{eq:Pcov} that $\Gamma_{\mathrm{post}, P} \geq \Gamma_{\mathrm{post}}$ in the sense of positive-definiteness, which seems to indicate that the reconstruction based on \eqref{eq:Bayesian_system}, i.e.,~$\widehat{x}_{\mathrm{post}}$, is more reliable than its projection-based counterpart $\widehat{x}_{\mathrm{post},P}$. This would indeed be true if the data were produced by the system matrix $J$ in accordance with \eqref{eq:Bayesian_system}. However, if the data $y$ also includes contributions from, say, changes in the coupling coefficients or the absorption of the RONI, then \eqref{eq:mean} and \eqref{eq:cov} are based on a false assumption on the generation of data and thus provide misleading information. The formulas \eqref{eq:Pmean} and \eqref{eq:Pcov} try to compensate for this mismodeling by considering the reduced model \eqref{eq:proj_eq}, where the purpose of $P$ is projecting away the contributions by the unknown nuisance parameters. However, \eqref{eq:Pmean} and \eqref{eq:Pcov} are also based on a false assumption if the projection $P$ does not accomplish its task perfectly.  In addition, it should be kept in mind that neither of the models \eqref{eq:Bayesian_system} and \eqref{eq:proj_eq} accounts for the nonlinearity in the generation of the data $y$, making the corresponding estimates indicative at best.
\end{remark}

\section{Computational model}
\label{sec:experimental}
Our computational model is based on the segmented head model for a neonate at combined gestational and chronological age of 41.7 weeks from the database described in~\cite{collins2021construction} and developed by the University College London (UCL) and the Centre for the Developing Brain at King's College London, with data from the Developing Human Connectome Project (dHCP). Modifications to the initial orientation and segmentation of the head model are as detailed in~\cite{Hakula2026,Hirvi_2023}. Each voxel has physical dimensions of 1 $\times$ 1 $\times$ 1\, mm$^3$. 

\subsection{Optical parameters}
The head model is segmented into five tissue types: combined scalp and skull (S\&S), two cerebrospinal fluid (CSF) types, gray matter (GM), and white matter (WM). In particular, as in \cite{Hirvi_2023}, the cerebrospinal fluid (CSF) is divided into two layers by isolating the subarachnoid layer in the sulci and ventricles. See the right panel of Figure~\ref{fig:meas_geom} for more details. 

The tissue-specific optical parameters adapted for the selected wavelength of 798\,nm are given in Table \ref{table:param}. The values were approximated with literature-based estimates provided at 800\,nm ~\cite{fukui2003,jonsson2018,Hirvi_2023}.

\begin{table}[t!]
\centering
\caption{Selected optical parameters for the considered tissue types. CSF-1 refers to the semidiffusive cerebrospinal fluid (CSF) in the subarachnoid region, whereas CSF-2 is the clearer CSF in the sulci and ventricles.}  
\begin{tabular}{lcccc}
 Tissue Type & $\mu_{a}$ [mm$^{\mathrm -1}$] & $\mu_{s}$ [mm$^{\mathrm -1}$] & $g$ & $\nu$\\ 
 \midrule
  Scalp \& Skull & 0.015 & 16 & 0.9 & 1.4 \\
  CSF-1 & 0.004 & 1.6 & 0.9 & 1.4 \\
  CSF-2 & 0.002 & 0.4 & 0.9 & 1.4 \\
  Gray Matter & 0.048 & 5 & 0.9 & 1.4 \\
  White Matter & 0.037 & 10 & 0.9 & 1.4 \\	 
 \bottomrule	
\end{tabular}
\label{table:param}
\end{table}

\subsection{Optode placement}
\begin{figure}[t!]
    \centering
    \includegraphics[width=\textwidth]{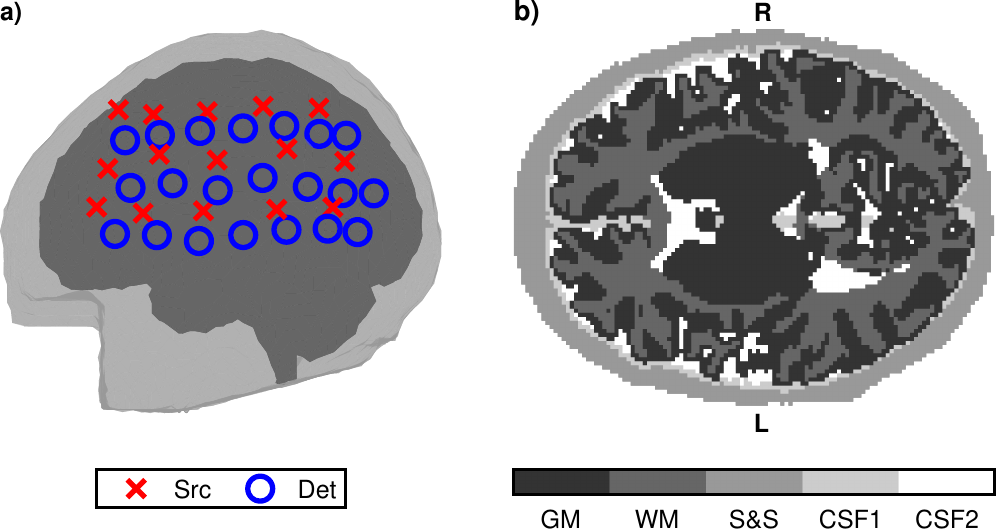}
    \caption{a) The 15 sources (red crosses) and 21 detectors (blue rings) on the studied neonate head. Face has been cut for anonymization. b) Axial slice of the head visualized with the segmented tissue types: GM = gray matter, WM = white matter, S\&S = combined scalp and skull, CSF1 = cerebrospinal fluid (CSF) in the subarachnoid region, CSF2 = CSF in the sulci and ventricles.}
    \label{fig:meas_geom}
\end{figure}

The employed optode configuration, consisting of 15 sources and 21 detectors placed over the left hemisphere of the neonate, is shown on the left panel in Figure~\ref{fig:meas_geom}. The shortest source--detector separation (SDS) is about 5\,mm, and we only consider optode pairs with an SDS below 40\,mm to avoid simulations with high stochastic noise. This results in a total of $m=210$ source--detector pairs and $l = 36$ coupling coefficients. We employed a collimated Gaussian beam source with a waist radius (half-width) of 1.25\,mm inserting intensity-modulated light at the modulation frequency of 100\,MHz. Each source launches photons in the direction of the surface normal in the source area, with the spatial launch probability following the Gaussian profile. The isotropic detection area was selected as the intersection of a sphere with radius of 1.82\,mm with the exterior head surface. The same settings were previously applied for an actual measurement probe with a similar structure in~\cite{maria2022}. Reflections at light input are neglected, but back-reflections are considered when a photon packet attempts to exit at the tissue--air boundary.

\subsection{Jacobians and the field-of-view}
\label{sec:FOV}
All computations were performed with the latest release (pre-compiled binaries for version 2025.10 Kilo-Kelvin) of the MCX software. The measurements were simulated according to \eqref{eq:weight}--\eqref{eq:amp_phase}, with billion photon packets for each source. The Jacobians for the log-amplitude and phase with respect to changes in the per-voxel absorption coefficients were computed in the ``rf replay'' mode of MCX; for details see~\cite[Section~2.2.1]{Hirvi_2023}. Importantly, in the ``replay'' mode (combined with mBLL), only one forward simulation per source is required to record the trajectories for the detected photon packets, and these trajectories can be ``replayed'' efficiently (and for different sets of the baseline absorption coefficients) to compute the Jacobians~\cite{YaoReplay,Hirvi_2023}.  

The difference data and Jacobians for each source were computed in parallel with one 16/32/80\,GB NVIDIA V100/A100/H100 GPU card randomly assigned for each array job on the 64-bit Linux-based Triton computing cluster provided by the Science-IT project of Aalto University School of Science. A single-source forward simulation for the $10^9$ photon packets executed in a couple of minutes.

The FOV was limited to the voxels for which the sensitivity is over 1\% of the maximum sensitivity in the brain for the corresponding measurement for any source--detector pair,~cf.~\cite{maria2022}.

\subsection{Prior and noise models}
For absorption perturbations, we assume the prior covariance structure
\begin{equation*}
\left( \Gamma_{x_{\rm total}}^0 \right)_{i,j} = \sigma^2\, \exp{ \left( -\frac{\| z_i - z_j \|_2^2}{2d^2} \right)}, \qquad i, j = 1, \dots, n_{\rm total},
\end{equation*}
where $\sigma$ is the pointwise standard deviation, $d$ is the spatial correlation length, and $z_i$ and $z_j$ are the coordinates of the voxels with indices $i$ and $j$ in the FOV. We select $\sigma$~=~0.003\,mm$^{-1}$ and $d$~=~3\,mm in all our numerical experiments; the former represents the expected contrast of perturbations in the absorption coefficient, whereas the latter reflects our {\em a priori} knowledge about the diameters of activated volumes in the brain. To promote fast computations via sparse structures, individual covariance values that are less than 0.01\% of $\sigma^2$ are set to zero, leading to
\begin{equation}
\label{eq:truncation}
\left( \Gamma_{x_{\rm total}} \right)_{i,j} = 
\left\{
\begin{array}{ll}
\left( \Gamma_{x_{\rm total}}^{0} \right)_{i,j} \quad &\mathrm{if\ } \left( \Gamma_{x_{\rm total}}^{0} \right)_{i,j} > (0.01 \, \sigma)^2\, , \\[1mm]
0 \quad  &\mathrm{otherwise}.
\end{array}
\right.
\end{equation}
Note that this covariance matrix $\Gamma_{x_{\rm total}} \in \R^{n_{\rm total} \times n_{\rm total}}$ corresponds to all voxels in the FOV, while the covariance matrices for restrictions to the voxels in the ROI and RONI, i.e.~$\Gamma_x \in \R^{n \times n}$ and $\Gamma_{\tilde{x}} \in \R^{\tilde{n} \times \tilde{n}}$, are obtained as the respective diagonal blocks of $\Gamma_{x_{\rm total}}$~\cite{kaipio2005book}. Moreover, observe that even though the truncation in \eqref{eq:truncation} may affect the positive definiteness and invertibility of~$\Gamma_{x}$, this induces no direct consequence on the usability of \eqref{eq:mean}--\eqref{eq:cov} and \eqref{eq:Pmean}--\eqref{eq:Pcov} because they do not involve operating with~$\Gamma_x^{-1}$. 

The components of the additive zero-mean Gaussian noise process are assumed to be independent, which leads to a diagonal noise covariance $\Gamma_e$. Recall that the first $m$ elements in $y \in \R^{2m}$ correspond to log-amplitude difference measurements and the last $m$ to phase difference measurements. We define the standard deviations for the two data types as
\[
\gamma_{\ln \! A} = 0.01\! \max_{1 \leq j \leq m} |y^0_{j}| \qquad \text{and} \qquad \gamma_{\varphi} = 0.01\! \! \! \! \max_{m+1 \leq j \leq 2m} |y^0_{j}|, 
\]
where $y^0 \in \R^{2 m}$ is the simulated difference data vector with no additive noise. Hence,
\[
\Gamma_e = \text{diag}\big( \gamma_{\ln \! A}^2 I, \gamma_{\varphi}^2 I \big),
\]
where $I \in \R^{m \times m}$ is an identity matrix.
Loosely speaking, this means that in our numerical experiment, the difference data contains 1\% of additive Gaussian noise. In particular, since the inversion does not directly suffer from stochastic noise related to the MC simulation because the Jacobian matrices are computed using the same photon paths as used for the simulation of data, adding artificial noise to the difference data is essential for avoiding an inverse crime.

\section{Numerical experiments}
\label{sec:numerics}
\begin{figure}[t!]
    \centering
    \includegraphics[width=\linewidth]{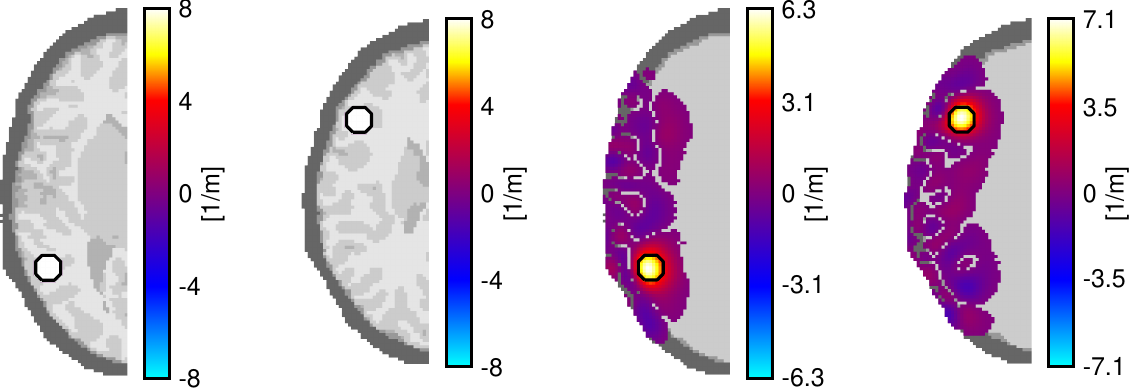} 
    \caption{Two left columns: cross-sections of target~1. Two right columns: the corresponding cross-sections of the reference reconstruction with $L^2$ error $0.185$ over the field-of-view.}
    \label{fig:Target1}
\end{figure}

\begin{figure}[t!]
    \centering
    \includegraphics[width=\linewidth]{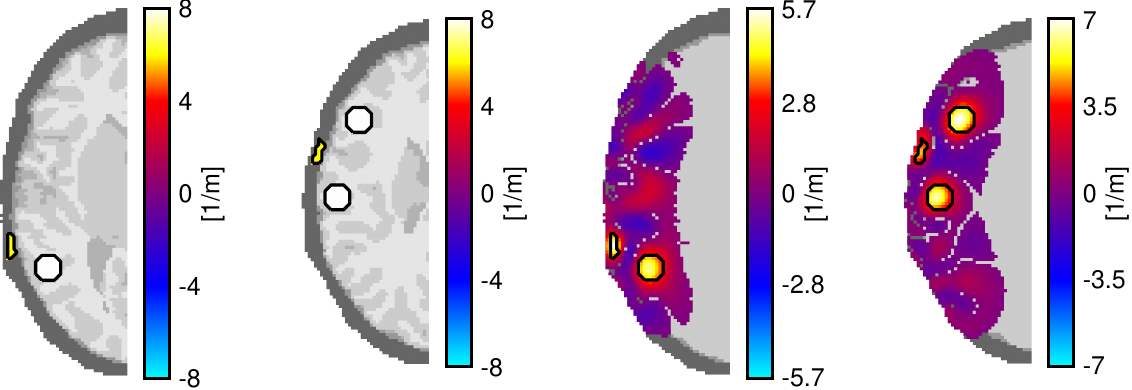}
    \caption{Two left columns: cross-sections of target~2. Two right columns: the corresponding cross-sections of the reference reconstruction with $L^2$ error $0.266$ over the field-of-view.}
    \label{fig:Target2}
\end{figure}

In this section, we present the results of our numerical experiments. We consider the two targets presented in the two leftmost columns of Figures~\ref{fig:Target1} and \ref{fig:Target2}, which correspond to two axial slices across the perturbed regions at different heights in the same head model. The first one carries two absorption perturbations in the brain with contrast 8\,m$^{-1}$ and none in the S\&S, whereas the second one has three absorption perturbations in the brain with contrast 8\,m$^{-1}$ and two in the S\&S with contrast 6\,m$^{-1}$. Each perturbation has a radius of 5\,mm. The perturbation sizes and contrasts were adapted from~\cite{Hakula2026}, and similar high-contrast brain perturbations were used in~\cite{heiskala2009probab, heiskala2009significance}. The corresponding standard reconstructions, i.e.~$\widehat{x}_{\mathrm{post}}$ from \eqref{eq:mean}, are visualized in the two rightmost columns of Figures~\ref{fig:Target1} and \ref{fig:Target2}, with the assumption that the ROI covers the whole FOV and there are no uncertainties in the coupling coefficients nor the baseline values for the optical parameters. Since we consider the same two targets in all following examples, the reconstructions in Figures~\ref{fig:Target1} and \ref{fig:Target2} serve as our benchmarks; their $L^2$ reconstruction errors over the FOV are $0.185$ and $0.266$, respectively. The leftmost figures visualize all tissue types in the background, but for the rest of the figures in this section, the tissue types inside the brain (GM, WM, and CSF-2) are visualized as one.

\subsection{Case 1: Mismodeled coupling coefficients}
\label{sec:case1}
We assume the optodes in the reference measurements are ideally coupled.
%
For the measurements simulated for the targets in Figures~\ref{fig:Target1} and \ref{fig:Target2}, we perturb the coupling coefficients at all optodes such that
\[
    A^{(s)},A^{(d)} \sim U(\delta_{A},1), \quad \varphi^{(s)},\varphi^{(d)} \sim U(0, \delta_{\varphi}),
\]
where $\delta_{A} = 0.9$ and $\delta_{\varphi} = \pi / 360$\, rad, and $U(a,b)$ denotes the uniform probability distribution on the interval $[a,b]$. The selected intervals were from~\cite{mozumder2013optode}. This leads to a mismatch in the coupling coefficients between the two measurements forming the difference data for both targets. The coupling Jacobian $J_c$ with respect to the log-amplitudes and phases of the coupling coefficients is formed as explained in Section~\ref{sec:ccs}, and the corresponding projection matrix $P$ onto the orthogonal complement of the range of $J_c$ is computed as indicated in~\eqref{eq:proj}. The ROI is set to coincide with the FOV.

\begin{figure}[!t]
    \centering
    \includegraphics[width=\linewidth]{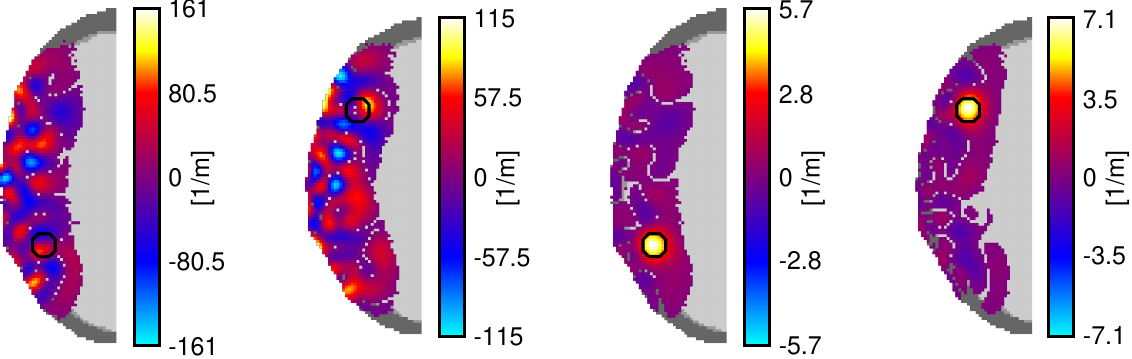}
    \caption{Case~1 with target~1. Two left columns: reconstruction cross-sections without a projection with $L^2$ error $9.940$ over the field-of-view (FOV). Two right columns: reconstruction cross-sections utilizing a projection with respect to the coupling coefficients with $L^2$ error $0.195$ over the FOV.}
    \label{fig:Proj_CC}
\end{figure}

\begin{figure}[!t]
    \centering
    \includegraphics[width=\linewidth]{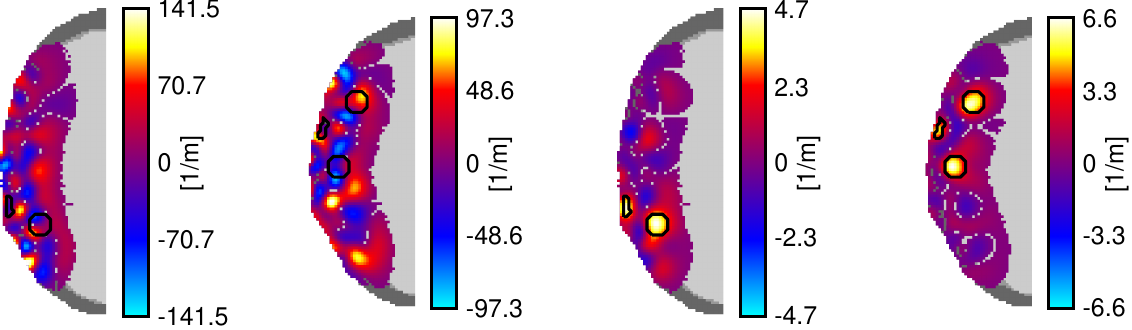}
    \caption{Case~1 with target~2. Two left columns: reconstruction cross-sections without a projection with $L^2$ error $8.562$ over the field-of-view (FOV). Two right columns: reconstruction cross-sections utilizing a projection with respect to the coupling coefficients with $L^2$ error $0.266$ over the FOV.}
    \label{fig:Proj_CC2}
\end{figure}

The resulting reconstructions for the two targets without and with $P$, i.e.~$\widehat{x}_{\mathrm{post}}$ from \eqref{eq:mean} and $\widehat{x}_{\mathrm{post},P}$ from \eqref{eq:Pmean}, are illustrated in the leftmost two and the rightmost two columns of Figures~\ref{fig:Proj_CC} and~\ref{fig:Proj_CC2}, respectively. In both cases, not accounting for the changes in the coupling coefficients results in completely uninformative reconstructions, whereas including the orthogonal projection in \eqref{eq:Pmean} yields reconstructions that are comparable in quality to the baseline reconstructions in Figures~\ref{fig:Target1} and \ref{fig:Target2}. This is verified by the corresponding $L^2$ reconstruction errors that are $0.195$ and $0.266$, respectively, i.e., practically the same as for the baseline reconstructions.

\subsection{Case~2: Uninteresting perturbations in the RONI}
\label{sec:case2}
In functional DOT studies, the focus is on detecting hemodynamic responses of activation in the brain or part of the brain. When the experiment consists of long time series of data and large 3D volumes of unknown voxels, restricting the reconstruction to the ROI leads to a significant dimension reduction. 

Next, we consider projecting away the effect of absorption perturbations outside the ROI (i.e.,~in the RONI) on the reconstruction in the ROI. For target~1 in Figure~\ref{fig:Target1}, the RONI is defined to be the frontal (i.e.,~top) half of the FOV, with the aim of demonstrating how projections can be employed to truncate the computational domain. For target~2 in Figure~\ref{fig:Target2}, we test whether a reconstruction in the brain can be computed as (or more) accurately by only considering voxels in the brain
and projecting away the influence of the perturbations in the S\&S, which we thus dub the RONI in this case.

\begin{figure}[!t]
    \centering
    \includegraphics[width=\linewidth]{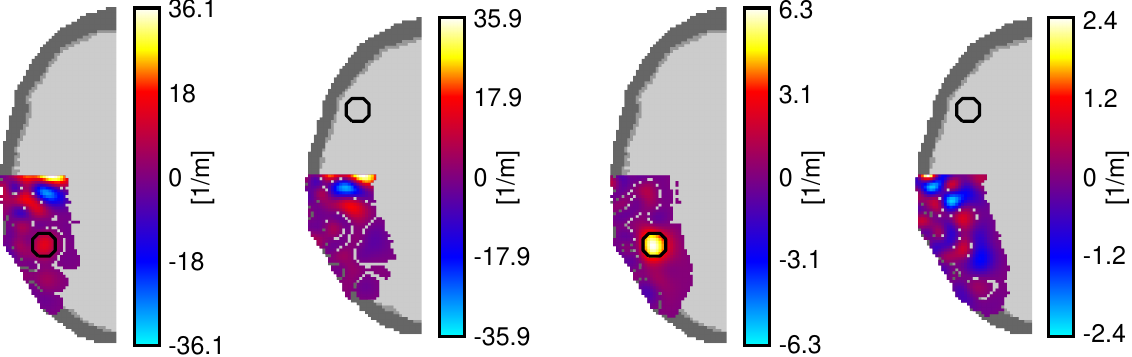}
    \caption{Case~2 with target~1 and the region-of-non-interest (RONI) being the frontal/top half of the field-of-view (FOV). Two left columns: reconstruction cross-sections without a projection with $L^2$ error $1.524$ over the ROI. Two right columns: reconstruction cross-sections utilizing a projection with respect to the RONI with $L^2$ error $0.183$ over the ROI.}
    \label{fig:Proj_RONI}
\end{figure}

\begin{figure}[!t]
    \centering
    \includegraphics[width=\linewidth]{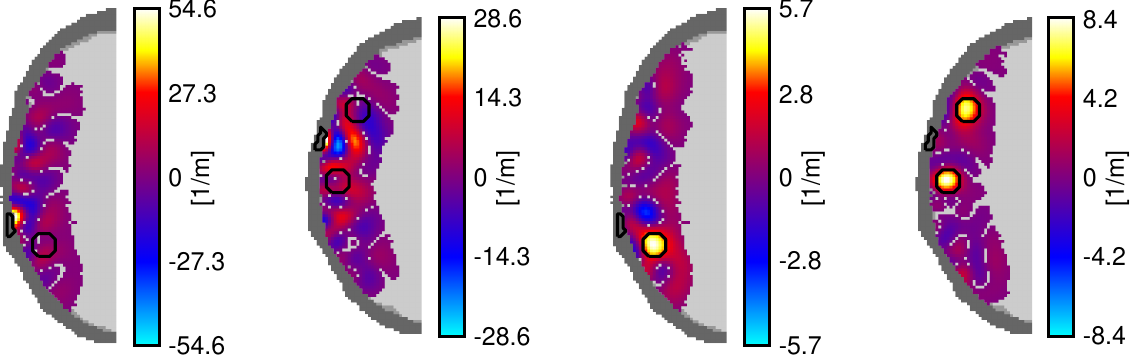}
    \caption{Case~2 with target~2 and the RONI being the S\&S. Two left columns: reconstruction cross-sections without a projection with $L^2$ error $1.220$ over the ROI. Two right columns: reconstruction cross-sections utilizing a projection with respect to the RONI with $L^2$ error $0.250$ over the ROI.}
    \label{fig:Proj_SS}
\end{figure}

For both targets, we build the orthogonal projection $P$ for reducing the effect of the absorption perturbations in the RONI by defining its nullspace to be the span of the eigenvectors of $\tilde{J} \Gamma_{\tilde{x}} \tilde{J}^{\rm \top}$ corresponding to the largest quarter (i.e.,~105) of the eigenvalues, cf.~Section~\ref{sec:high}. The reconstructions without and with $P$,~i.e.~$\widehat{x}_{\mathrm{post}}$ and $\widehat{x}_{\mathrm{post},P}$, are visualized, respectively, in the two leftmost and two rightmost columns of Figures~\ref{fig:Proj_RONI} and~\ref{fig:Proj_SS}, with the former corresponding to target~1 and the latter to target~2. As the naive reconstruction~$\widehat{x}_{\mathrm{post}}$, defined by \eqref{eq:mean}, tries to explain all difference data by absorption changes in the ROI, it exhibits in both cases strong artifacts along the boundary between the ROI and the RONI, rendering the reconstructions on the left in Figures~\ref{fig:Proj_RONI} and \ref{fig:Proj_SS} almost useless. On the other hand, since $\widehat{x}_{\mathrm{post}, P}$ defined by \eqref{eq:Pmean} compensates for the truncation of the reconstruction domain with the help of the projection $P$, the reconstructions of the ROI in the right columns of Figures~\ref{fig:Proj_RONI} and \ref{fig:Proj_SS} are almost as good as in the reference reconstructions of Figures~\ref{fig:Target1} and \ref{fig:Target2}.

These conclusions are supported by the $L^2$ reconstruction errors over the ROI, which are $0.183$ and $0.250$, respectively, for the reconstructions on the right in Figures~\ref{fig:Proj_RONI} and \ref{fig:Proj_SS}. As the respective numbers for the reference reconstructions on the right in Figures~\ref{fig:Target1} and \ref{fig:Target2} are $0.136$ and $0.244$, there is a moderate increase in the reconstruction error over the ROI compared to reference reconstruction for target~1 (RONI the frontal half of the FOV) and almost no increase for target~2 (RONI the S\&S layer).

\subsection{Case~3: Misspecified baseline parameters}
\label{sec:case3}
As an example of misspecified baseline optical parameters, we consider our estimates for the absorption coefficients in both GM and WM being off by a considerable margin in the computation of the Jacobian. More precisely, the assumed absorption levels when computing the Jacobians are as listed in Table~\ref{table:param}, while the absorption coefficients of the GM and WM with which the difference data is simulated were selected from~\cite{maria2022} as 0.014\,mm$^{-1}$ and 0.0032\,mm$^{-1}$, respectively, i.e., these parameters are applied to both the reference measurements and the measurements corresponding to the perturbations in Figures~\ref{fig:Target1} and~\ref{fig:Target2}. This corresponds to the realistic situation where the actual measured difference signal corresponds to different optical parameters than the ones used to simulate the Jacobians. Such errors in the specifications of the optical parameters for the GM and WM are indeed possible based on the values found in literature; see~\cite{farina2015invivo,maria2022,mozumder2024diffuse}. However, it should also be noted that the reconstruction process appeared quite robust to the baseline absorption levels in the considered simulated examples, and hence visually demonstrating an improvement in the localization of the absorption changes via the projection approach requires relatively significant discrepancy in the baseline values. Apart from the absorption levels of the GM and WM, the other optical parameters are the same for the simulation of data and formation of the absorption Jacobian $J$ for the ROI, which is defined to be the FOV in this example. Note that the FOV here is the same as in all other examples since it is defined with the Jacobian for the assumed (not the true) absorption levels.

\begin{figure}[!t]
    \centering
    \includegraphics[width=\linewidth]{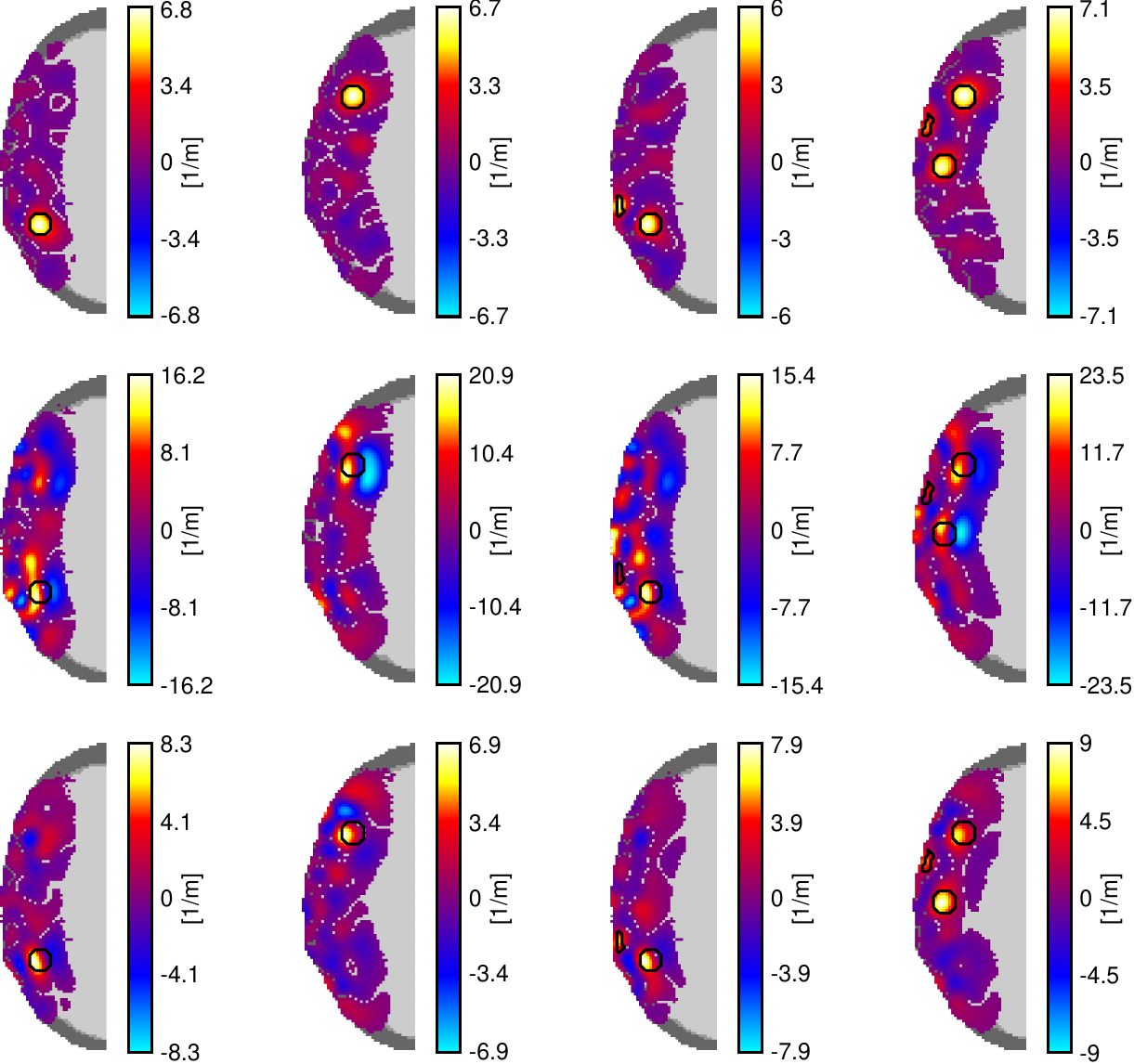}
    \caption{Case~3 with target~1 (two left columns) and target~2 (two right columns). Top row: reconstruction cross-sections using the accurate absorption levels for the GM and WM with $L^2$ errors $0.211$ (left) and $0.273$ (right) over the field-of-view (FOV). Middle row: reconstruction cross-sections using misspecified absorption levels for the GM and WM without a projection with $L^2$ errors $1.172$ (left) and $1.232$ (right) over the FOV. Bottom row: reconstruction cross-sections using the misspecified absorption levels {\rm and} a projection with respect to the absorption levels of the GM and WM with $L^2$ errors $0.387$ (left) and $0.399$ (right) over the FOV.}
    \label{fig:Proj_Baseline}
\end{figure}

Adopting the notation of Section~\ref{sec:uncertain} and specifically that of \eqref{eq:sum_matrix}, the first `difference matrix' $\hat{J}_{1,\delta_1}$ is simulated with the absorption coefficient 0.031\,mm$^{-1}$ for the GM (i.e., with $\delta_1= -0.017$\,mm$^{-1}$) and the second one $\hat{J}_{2,\delta_2}$ with the absorption coefficient 0.024\,mm$^{-1}$ for the WM (i.e., with $\delta_2=  -0.013$\,mm$^{-1}$). In particular, note that the employed perturbed values for the absorption levels of the GM and WM are much closer to the assumed absorption coefficients than to the true ones used for simulating the data. To form the projection $P$, we consider the matrix \eqref{eq:sum_matrix} with $q=2$ and $A = \Gamma_x$ and use its eigenvectors associated to the largest 105 eigenvalues, i.e., the first fourth of the eigenvalues, as the basis for the nullspace of $P$, as in case 2.

The results for the two targets are shown in Figure~\ref{fig:Proj_Baseline}, which is organized in a slightly modified manner compared to the figures presenting the findings of cases~1 and~2. Here, the two leftmost columns correspond to target~1 and the two rightmost columns to target~2. The top row visualizes the reconstructions obtained using the accurate baseline parameters, i.e., the reconstructions obtained with the matrix $J_{\Delta}$ from \eqref{eq:uncertain} replacing $J$ in \eqref{eq:mean}. Thus, the first row works as our reference in this case; the differences compared to Figures~\ref{fig:Target1}--\ref{fig:Target2} stem from the different baseline optical parameters. The middle row visualizes the reconstructions obtained with the mismodeled absorption levels of the GM and WM, and the bottom row corresponds to the reconstruction with the projection $P$. If $P$ is not involved in the reconstruction process, the discrepancy in the absorption levels of the GM and WM between the simulation of the measurements and the computation of the Jacobian~$J$ used in \eqref{eq:mean} clearly affects the obtained reconstructions for both targets. The perturbations in the brain are located with quite a good accuracy, but the perturbations near the surface in target~2 are not visible. The contrast of the reconstructed perturbations is also too high and there are significant undulations in the background. The reconstructions with $P$, presented in the bottom row, also suffer from background noise when compared to the reference reconstructions on the top row, but the effect is milder than in the middle row. However, the perturbations are better, yet not perfectly located, and their contrast is more reasonable.

The $L^2$ errors confirm that the adoption of $P$ improves the localization of the absorption perturbations, though not as significantly as in cases~1 and 2. In Figure~\ref{fig:Proj_Baseline}, the top-row reference reconstructions yield errors of $0.211$ (target~1) and $0.273$ (target~2), the middle-row naive reconstructions yield $1.172$ and $1.232$, and the bottom-row reconstructions utilizing $P$ yield $0.387$ and $0.399$, respectively.

\subsection{Case~4: Multiple uncertainties}
As an example on using projections to reduce the effects of several uncertainties in the measurement model, we combine cases~1 and 2. That is, we induce discrepancy in the coupling coefficients between the reference measurements and the measurements simulated for the targets in Figures~\ref{fig:Target1} and \ref{fig:Target2} and form the reconstruction in the ROI only, even though there are also absorption perturbations in the RONI. We manipulate the coupling coefficients for the simulations of measurements with targets~1 and 2 as explained in Section~\ref{sec:case1}, and the specifications of the ROI and RONI for the two targets are the same as in case~2 described in Section~\ref{sec:case2}.

For each target, the orthonormal basis vectors for the nullspace of the orthogonal projection for reducing the effect of RONI is formed as explained in Section~\ref{sec:case2} with $105$ eigendirections, and the vectors are then stacked as the columns of a matrix $\tilde{V}$. Following the instructions in Section~\ref{sec:combo}, the matrix $\mathcal{J} = [ J_c, \tilde{V}]$, with $J_c$ being the coupling Jacobian from Section~\ref{sec:ccs}, is then used in place of mere $J_c$ in \eqref{eq:proj} when forming the orthonormal projection $P$ to be used in the reconstruction formulas of Section~\ref{sec:background}.

\begin{figure}
    \centering
    \includegraphics[width=\linewidth]{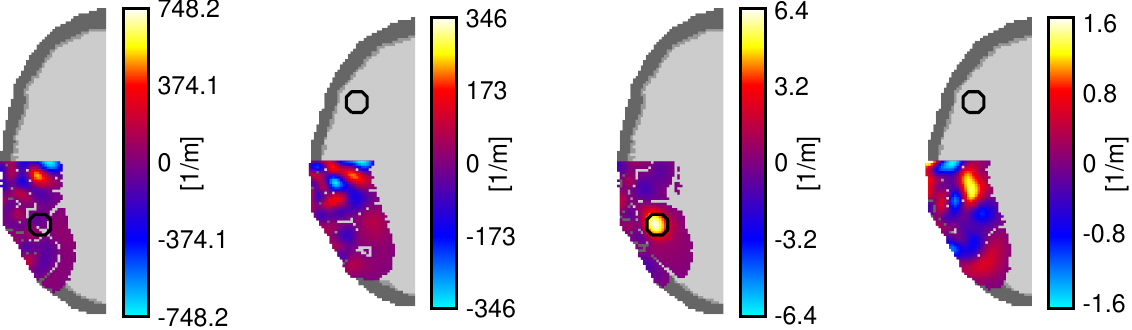}
    \caption{Case~4 with target~1. Two left columns: reconstruction cross-sections without a projection with $L^2$ error $23.730$ over the ROI. Right: reconstruction cross-sections utilizing a projection with respect to the coupling coefficients and the RONI with $L^2$ error $0.158$ over the ROI.}
    \label{fig:Proj_multi1}
\end{figure}

\begin{figure}
    \centering
    \includegraphics[width=\linewidth]{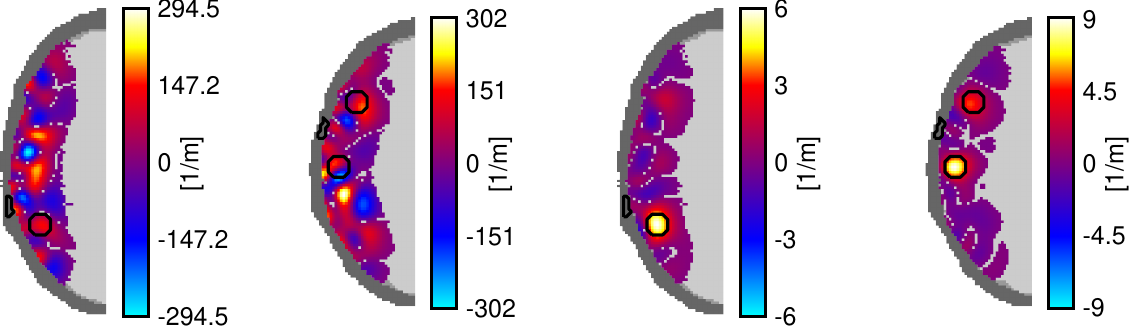}
    \caption{Case~4 with target~2. Two left columns: reconstruction cross-sections without a projection with $L^2$ error $16.290$ over the ROI. Right: reconstruction cross-sections utilizing a projection with respect to the coupling coefficients and the RONI with $L^2$ error $0.238$ over the ROI.}
    \label{fig:Proj_multi2}
\end{figure}

The reconstructions for the two targets without and with $P$ are presented in the two leftmost and two rightmost columns of Figures~\ref{fig:Proj_multi1} and \ref{fig:Proj_multi2}. In both figures, the reconstructions not employing $P$ are even less informative than the corresponding ones in Section~\ref{sec:case1} and \ref{sec:case2}, where the effect of alterations in the coupling coefficients and uninteresting absorption changes in the RONI were separately considered. The reconstructions accounting for the mismodeling via the use of $P$ are arguably somewhat worse than the corresponding ones in Sections~\ref{sec:case1} and \ref{sec:case2}, but they still exhibit a clear improvement compared to the naive reconstructions. The $L^2$ errors over the respective ROIs for the reconstructions on the right in Figures~\ref{fig:Proj_multi1} and \ref{fig:Proj_multi2} are $0.158$ and $0.238$, respectively. These numbers are approximately the same as the corresponding $L^2$ errors over the appropriate ROI's for the reconstructions on the right in Figures~\ref{fig:Target1} and \ref{fig:Target2},  Figures~\ref{fig:Proj_CC} and \ref{fig:Proj_CC2}, and Figures~\ref{fig:Proj_RONI} and \ref{fig:Proj_SS}, respectively.

\section{Concluding remarks}
\label{sec:conclusion}
This work applied the projection technique introduced in~\cite{Calvetti2025,Jaaskelainen25,Jaaskelainen26} for handling unknown nuisance parameters in inverse problems to the framework of brain imaging by DOT. The presented numerical experiments demonstrated that projecting the linearized measurement model of DOT onto the orthogonal complement of a subspace that is expected to be affected the most by the nuisance parameters can considerably reduce artifacts caused by mismodeled coupling coefficients, truncation of the computational domain, or a mispecified tissue-wise baseline absorption level.

For mismodeled coupling coefficients or truncation of the computational domain, the aforementioned subspace that serves as the nullspace of the utilized projection was formed by considering the properties of the Jacobian matrix with respect to the considered nuisance parameter,~i.e.,~the coupling coefficients or the absorption change in the RONI. A similar idea was used in~\cite{Jaaskelainen25,Jaaskelainen26} to address domain truncation and misspecified contact resistances in EIT, and our numerical results are well-aligned with those in~\cite{Jaaskelainen25,Jaaskelainen26}. However, the forward model in ~\cite{Jaaskelainen25,Jaaskelainen26} was defined by an elliptic partial differential equation, while our numerical computations for DOT were based on MC simulations of the RTE. 

According to our knowledge, our technique for building the projection for reducing the effect of a misspecified tissue-wise baseline absorption level has not been employed previously for any imaging modality. The nullspace for the projection was constructed by considering eigenvectors of a (weighted) difference of two Jacobian matrices with respect to the voxelized absorption in the ROI evaluated at two different absorption levels for the considered tissue type. It can be argued that such a difference, in fact, approximates the derivative of the absorption Jacobian with respect to the baseline absorption value in question. Investigating whether the employment of such a second derivative in projecting away the effect of misspecified tissue-wise optical parameters (including the scattering coefficient) improves the numerical results provides a possible avenue for future research. The same applies to introducing explicit rules for choosing the dimension of the nullspace for nuisance projections.

\section*{Acknowledgements}
We acknowledge the computational resources provided by the Aalto Science-IT project. We utilized the head model data made available from the Developing Human Connectome Project funded by the European Research Council under the European Union's Seventh Framework Programme (FP/2007-2013)/ERC Grant Agreement no.\ (319456).

\appendix

\bibliographystyle{siamplain}
\bibliography{biblio.bib}

\end{document}